\documentclass[12pt]{article}
\usepackage{amsmath}
\usepackage{amsfonts}
\usepackage{graphicx}
\usepackage{amssymb}
\usepackage{theorem}

\usepackage[OT2,OT1]{fontenc}
\newcommand\cyr{%
	\renewcommand\rmdefault{wncyr}%
	\renewcommand\sfdefault{wncyss}%
	\renewcommand\encodingdefault{OT2}%
	\normalfont
	\selectfont}
\DeclareTextFontCommand{\textcyr}{\cyr}

\newtheorem{theorem}{Theorem}

\newtheorem{corollary}{Corollary}

\newtheorem{proposition}{Proposition}
\newtheorem{remark}{Remark}

\usepackage[russian, english]{babel}
\usepackage{centernot}
 
\usepackage[unicode,colorlinks,plainpages=false]{hyperref}
\usepackage{comment}
 
\newcommand{\openbox}{$\begin{array}{c}
		\hspace*{-0.55em}\sqcap \hspace*{-0.60em}\\[-0.4em] \hline
		\multicolumn{1}{c}{\hspace*{-0.60em}}\\[-0.8em]
	\end{array}$}

\begin{document}
			
			\centerline{\bf  Right regular triples of semigroups}

			\bigskip
			
			\centerline{Csaba T\'oth}
			
			\bigskip
			
			\centerline{Department of Algebra, Institute of Mathematics}
			\centerline{Budapest University of Technology and Economics}
			\centerline{1521 Budapest, P.O. Box 91, Hungary}
			\centerline{e-mail: tcsaba94@gmail.com}

			\bigskip
			
			\begin{abstract}
			Let ${\cal M}(S; \Lambda; P)$ denote a Rees $I\times \Lambda$ matrix semigroup without zero over a semigroup $S$, where $I$ is a singleton. If $\theta _S$ denotes the kernel of the right regular representation of a semigroup $S$, then a triple $A, B, C$ of semigroups is said to be right regular, if there are mappings $A\stackrel{P}{\longleftarrow}B$ and $B\stackrel{P'}{\longrightarrow}C$ such that
			${\cal M}(A; B; P)/\theta_{{\cal M}(A; B; P)}\cong {\cal M}(C; B; P')$. In this paper we examine right regular triples of semigroups.
			\end{abstract}
			
			\bigskip
			
			{Keywords}: semigroup, congruence, Rees matrix semigroup
			
			\medskip
			
			{2010 Mathematics Subject Classification}: 20M10, 20M30.
			
			\bigskip

\section{Introduction and motivation}\label{intr}

The notion of right regular triples of semigroups is defined in \cite{NagyToth}, where a special type of Rees matrix semigroups without zero over semigroups are examined. A triple $A, B, C$ of semigroups is said to be right regular, if there are mappings \[A\stackrel{P}{\longleftarrow}B\stackrel{P'}{\longrightarrow}C\] such that the factor semigroup ${\cal M}(A; B; P)/\theta _{{\cal M}(A; B; P)}$ is isomorphic to the semigroup ${\cal M}(C; B; P')$, where $\theta _{{\cal M}(A; B; P)}$ is the kernel of the right regular representation of the semigroup ${\cal M}(A; B; P)$. In \cite{NagyToth} it is proved that if $A, B, C$ are semigroups such that $A/\theta _A \cong B$ and $B/\theta _B\cong C$, then the triple $A, B, C$ is right regular. There is also an example given for a right regular triple $A, B, C$ of semigroups such that none of the conditions $A/\theta _A\cong B$ and $B/\theta _B\cong C$ are fulfilled.  These results motivate us to investigate right regular triples of semigroups. In this paper we examine the connection between the structure of semigroups belonging to a right regular triples of semigroups, and present quite a few examples of right regular triples of semigroups.

\section{Preliminaries}
By a semigroup we mean a multiplicative semigroup, that is, a nonempty set endowed with an associative multiplication.
	
A nonempty subset $L$ of a semigroup $S$ is called a \emph{left ideal} of $S$ if $SL\subseteq L$. The concept of a right ideal of a semigroup is defined analogously. A semigroup $S$ is said to be \emph{left (resp., right) simple} if $S$ itself is the only left (resp., right) ideal of $S$. A semigroup $S$ is left (resp., right) simple if and only if $Sa=S$ (resp., $aS=S$) for every $a\in S$.
	
	A semigroup $S$ is called \emph{left cancellative} if $xa=xb$ implies $a=b$ for every $x, a, b\in S$. A left cancellative and right simple semigroup is called a right group. A semigroup satisfying the identity $ab=b$ is called a \emph{right zero semigroup}. By \cite[Theorem 1.27.]{Clifford1:sg-1}, a semigroup is a right group if and only if it is a direct product of a group and a right zero semigroup.
 
In \cite[Theorem 1]{Cohn:sg-1}, it is shown that a semigroup $S$ is embedded in an idempotent-free left simple semigroup if and only if $S$ is idempotent-free and satisfies the condition: for all $a, b, x, y\in S$, $xa=xb$ implies $ya=yb$.
Using the terminology of \cite{Nagy:sg-8}, a semigroup $S$ satisfying this last condition is called a \emph{left equalizer simple semigroup}. In other words, a semigroup $S$ is left equalizer simple if, for arbitrary elements $a, b\in S$, the assumption that $xa=xb$ is satisfied for some $x\in S$ implies that $ya=yb$ is satisfied for all $y\in S$. By \cite[Theorem 2.1]{Nagy:sg-8}, a semigroup $S$ is left equalizer simple if and only if the factor semigroup $S/\theta _S$ is left cancellative.
	
A nonempty subset $I$ of a semigroup $S$ is called an \emph{ideal} of $S$ if $I$ is a left ideal and a right ideal of $S$. A semigroup $S$ is called \emph{simple} if $S$ itself is the only ideal of $S$. By \cite[Lemma 2.2.]{Clifford1:sg-1}, a semigroup $S$ is simple if and only if $SaS=S$ for every $a\in S$.
	
	Let $S$ be a semigroup and $I$ be an ideal of $S$. We say that the homomorphism $\varphi: S \mapsto I$ is a \emph{retract homomorphism} \cite[Definition 1.44]{Nagybook:sg-5}, if it leaves the elements of $I$ fixed. In this case, $I$ is called a \emph{retract ideal} of $S$, and $S$ is a \emph{retract extension} of $I$ by the Rees factor semigroup $S/I$.

A transformation $\varrho$ of a semigroup $S$ is called a \emph{right translation} of $S$ if $(xy)\varrho=x(y\varrho)$ is satisfied for every $x, y\in S$. For an arbitrary element $a$ of a semigroup $S$, $\varrho _a: x\mapsto xa$ ($x\in S$) is a right translation of $S$ which is called an \emph{inner right translation} of $S$ corresponding to the element $a$. For an arbitrary semigroup $S$, the mapping $\Phi _S: a\mapsto \varrho _a$ is a homomorphism of $S$ into the semigroup of all right translations of $S$. The homomorphism $\Phi _S$ is called the \emph{right regular representation} of $S$. For an arbitrary semigroup $S$, let $\theta _S$ denote the kernel of $\Phi_S$. It is clear that $(a,b)\in \theta _S$ for elements $a, b\in S$ if and only if $xa=xb$ for all $x\in S$. A semigroup $S$ is called \emph{left reductive} if $\theta _S$ is the identity relation on $S$. Thus $\theta _S$ is faithful if and only if $S$ is left reductive. The congruence $\theta _S$ plays an important role in the investigation of the structure of the semigroup $S$.
 In \cite{Chrislock:sg-1}, the author characterizes semigroups $S$ for which the factor semigroup $S/\theta _S$ is a group. In \cite{Chrislock:sg-2}, semigroups $S$ are characterized for which the factor semigroup $S/\theta _S$ is a right group. In \cite[Theorem 2]{NagyKolibiar:sg-19}, a construction is given which shows that every semigroup $S$ can be obtained from the factor semigroup
$S/\theta _S$ by using this construction. In \cite{NagyTothprob}, the authors study the probability that two elements which are selected at random with replacement from a finite semigroup have the same right matrix.
	
If $S$ is a semigroup, $I$ and $\Lambda$ are nonempty sets, and $P$ is a $\Lambda \times I$ matrix with entries $P(\lambda , i)$, then the set ${\mathcal M}(S; I, \Lambda ;P)$ of all triples $(i, s, \lambda )\in I\times S\times \Lambda$ is a semigroup under the multiplication $(i, s, \lambda )(j, t, \mu )=(i, sP(\lambda , j)t, \mu )$. According to the terminology in \cite[\textsection 3.1]{Clifford1:sg-1}, this semigroup is called a \emph{Rees $I\times \Lambda$ matrix semigroup without zero over the semigroup $S$ with $\Lambda \times I$ sandwich matrix $P$}. In \cite{NagyToth},
Rees matrix semigroups ${\cal M}(S; I, \Lambda ;P)$ without zero over semigroups $S$ satisfying $|I|=1$ are in the focus. In our present paper we also use such type of Rees matrix semigroups, which will be denoted by ${\cal M}(S;\Lambda ;P)$. In this case the matrix $P$ can be considered as a mapping of $\Lambda$ into $S$, and so the entries of $P$ will be denoted by $P(\lambda)$.
If the element of $I$ is denoted by $1$, then the element $(1, s, \lambda )$ of ${\cal M}(S;\Lambda ;P)$ can be considered in the form
$(s, \lambda )$; the operation on ${\cal M}(S;\Lambda ;P)$ is
$(s, \lambda )(t, \mu )=(sP(\lambda)t, \mu)$.

	For notations and notions not defined but used in this paper, we refer the reader to books \cite{Clifford1:sg-1}, \cite{Howie:sg-11}, and \cite{Nagybook:sg-5}.

	\section{Results }

	\begin{theorem}\label{proprightsimple}If $A, B, C$ is a right regular triple of semigroups such that $A$ is right simple, then $C$ is also right simple.
	\end{theorem}
	
	\textit{Proof.} Assume that $A, B, C$ is a right regular triple of semigroups. Then there are mappings $P: B\mapsto A$ and $P': B\mapsto C$ such that
		\[{\cal M}(A; B; P)/\theta _{{\cal M}(A; B; P)} \cong {\cal M}(C; B; P').\]
		Assume that $A$ is right simple. Let $(a_1, b_1), (a_2, b_2)\in {\cal M}(A; B; P)$ be arbitrary elements. Since $A$ is right simple, we have $aP(b_1)A=A$, and so there is an element $\xi \in A$ such that $a_1P(b_1)\xi=a_2$ and \[(a_1, b_1)(\xi , b_2)=(a_2, b_2).\] Hence the Rees matrix semigroup ${\cal M}(A; B; P)$ is right simple. As every homomorphic image of a right simple semigroup is right simple, the Rees matrix semigroup ${\cal M}(C; B; P')$ is right simple. Let
		$c, \eta \in C$ be an arbitrary elements.
		Then, for any $b \in B$,
		$(c, b ){\cal M}(C;B;P')={\cal M}(C;B;P')$, and so
		\[(c, b )(u, v )=(\eta, b)\] for some $(u, v )\in {\cal M}(C; B; P')$. Hence $cP'(b)u=\eta$. Thus $cC=C$ for every $c\in C$. Then $C$ is right simple. \hfill\openbox

	\begin{theorem}\label{proprightgroup} If $A, B, C$ is a right regular triple of semigroups such that $A$ is a right group, then $C$ is also a right group.
	\end{theorem}
	
	\textit{Proof.} Assume that $A, B, C$ is a right regular triple of semigroups. Then there are mappings $P: B\mapsto A$ and $P': B\mapsto C$ such that
		\[{\cal M}(A; B; P)/\theta _{{\cal M}(A; B; P)} \cong {\cal M}(C; B; P').\]
		Assume that $A$ is a right group, that is, right simple and left cancellative. By the proof of Theorem ~\ref{proprightsimple}, the semigroups ${\cal M}(A; B; P)$ and $C$ are right simple.
		Let $(a, b), (a_1, b_1), (a_2, b_2)\in {\cal M}(A; B ; P)$ be arbitrary elements with
		\[(a, b)(a_1, b_1)=(a, b)(a_2, b_2).\]
		Then \[(aP(b)a_1, b_1)=(aP(b)a_2, b_2),\] that is,
		\[aP(b)a_1=aP(b)a_2\quad \hbox{and}\quad b_1=b_2.\]
		As $A$ is left cancellative, we get $a_1=a_2$, and so \[(a_1, b_1)=(a_2, b_2).\] Hence the semigroup ${\cal M}(A; B ; P)$ is left cancellative. As ${\cal M}(A; B ; P)$ is also right simple, it is a right group. From the left cancellativity of ${\cal M}(A; B ; P)$ it follows that $\theta _{{\cal M}(A; B; P)}=\iota _{{\cal M}(A; B; P)}$. Thus the semigroup ${\cal M}(C; B ; P')$ is left cancellative.
		Assume $xc_1=xc_2$ for elements $x, c_1, c_2\in C$. Let $b\in B$ be arbitrary. As $C$ is right simple, there are elements $u, v\in C$ such that $P(b)u=c_1$ and $P(b)v=c_2$.
		Thus \[xP(b)u=xP(b)v.\] Then, for an arbitrary $b'\in B$,
		\[(x, b)(u, b')=(x, b)(v, b')\]
		is satisfied in ${\cal M}(C; B ; P)$.
		As ${\cal M}(C; B ; P)$ is left cancellative, we get $u=v$, from which it follows that $c_1=c_2$. Hence $C$ is left cancellative. By the above, $C$ is right simple. Consequently $C$ is a right group.
	\hfill\openbox
	
	\begin{theorem}\label{propsimple}If $A, B, C$ is a right regular triple of semigroups such that $A$ is simple, then $C$ is also simple.
	\end{theorem}
	
	\textit{Proof.} Assume that $A, B, C$ is a right regular triple of semigroups. Then there are mappings $P: B\mapsto A$ and $P': B\mapsto C$ such that
		\[{\cal M}(A; B; P)/\theta _{{\cal M}(A; B; P)} \cong {\cal M}(C; B; P').\]
		Assume that $A$ is simple. Let $(a, b), (u, v)\in {\cal M}(A; B; P)$ and $z\in B$ be an arbitrary elements. Then $AP(z)aP(b)A=A$ implies that there are elements
		$\xi, \eta \in A$ such that $\xi P(z)aP(b)\eta =u$ and so \[(\xi , z)(a, b)(\eta , v)=(u, v).\] Hence the Rees matrix semigroup ${\cal M}(A; B; P)$ is simple. As every homomorphic image of a simple semigroup is simple, the Rees matrix semigroup ${\cal M}(C; B; P')$ is simple.
		
		Let $c_1, c_2\in C$ and $b_1, b_2\in B$ be arbitrary elements.
		Then \[{\cal M}(C; B; P')(c_1, b_1){\cal M}(C; B; P')={\cal M}(C; B; P'),\] and so there are elements $(x, \xi ), (y, \eta )\in {\cal M}(C; B; P')$ such that
		\[(xP(\xi )c_1P(b_1)y, \eta )=(x, \xi )(c_1, b_1)(y, \eta )=(c_2, b_2).\]
		Hence
		\[xP(\xi )c_1P(b_1)y=c_2.\] Thus \[Cc_1C=C\] for every $c_1\in C$. Then $C$ is simple. \hfill\openbox
	
\medskip

The next proposition is used in the proof of Theorem~\ref{leftcanc1}.

\begin{proposition}\label{lefteq} Let $A$ be a semigroup, $\Lambda$ be an arbitrary nonempty set and $P:\Lambda \mapsto A$ is an arbitrary mapping. If $A$ is left equalizer simple, then the Rees matrix semigroup ${\cal M}(A; \Lambda; P)$ is also left equalizer simple.
	\end{proposition}
	
	\textit{Proof.} Suppose that $A$ is a left equalizer simple semigroup, $\Lambda$ is a nonempty set and $P:\Lambda\mapsto A$ is a mapping. Take $(a_1,b_1),(a_2,b_2),(a,b) \in {\cal M}(A; \Lambda; P)$. Suppose that
		\[ (a,b)(a_1,b_1) = (a,b)(a_2,b_2). \]
		This means that
		\[ (aP(b)a_1,b_1) = (aP(b)a_2,b_2) \quad \iff \quad aP(b)a_1 = aP(b)a_2 \ \textrm{and} \ b_1 = b_2.  \]
		Since $A$ is left equalizer simple we have that, for all $x \in A$ and $y \in \Lambda:$
		\[\ xP(y)a_1 = xP(y)a_2,\]
		hence,
		\[ (x,y)(a_1,b_1) = (x,y)(a_2,b_2). \]
		Thus, ${\cal M}(A; \Lambda; P)$ is a left equalizer simple semigroup. \hfill\openbox

	\begin{theorem}\label{leftcanc1}Let $A, B, C$ be a right regular triple of semigroups
		such that $P': B \mapsto C $ is surjective. If $A$ is left equalizer simple, then $C$ is left cancellative.
	\end{theorem}
	
	\textit{Proof.} Assume that $A, B, C$ is a right regular triple of semigroups. Then there are mappings $P: B\mapsto A$ and $P': B\mapsto C$ such that
		\[{\cal M}(A; B; P)/\theta _{{\cal M}(A; B; P)}\cong {\cal M}(C; B; P').\]
		From Proposition~\ref{lefteq}, we have that ${\cal M}(A; B; P)$ is a left equalizer simple semigroup, and hence
		$ {\cal M}(C; B; P')$
		is left cancellative by \cite[Theorem 2.1]{Nagy:sg-8}. \\
		Now, take $x, c_1, c_2 \in C$ such that $xc_1 = xc_2$. Since $P'$ is surjective,   there exists $b \in B$ such that $P'(b) = x$. Then $P'(b)c_1 = P'(b)c_2$. Let $c\in C$ be arbitrary, then
		\[ (c,b)(c_1,b) = (cP'(b)c_1,b) = (cP'(b)c_2,b) = (c,b)(c_2,b). \]
		Since ${\cal M}(C; B; P')$ is left cancellative, $(c_1,b) = (c_2,b)$, hence $c_1 = c_2$. Thus $C$ is left cancellative. \hfill\openbox

	\medskip
	
	\begin{theorem}\label{leftcanc2} Let $A, B, C$ be a right regular triple of semigroups such that $C$ is left commutative. If $A$ is left equalizer simple, then $C$ is left cancellative.
	\end{theorem}
	
	\textit{Proof.} From the proof of Theorem \ref{leftcanc1}, we know that ${\cal M}(C; B; P')$ is left cancellative. Again, take $x,c_1,c_2 \in C$ such that $xc_1 = xc_2$. Then for arbitrary $b \in B$,
		\[P'(b)xc_1 = P'(b)xc_2.\]
		Since $C$ is left commutative,
		\[xP'(b)c_1 = xP'(b)c_2,\]
		and then
		\[ (x,b)(c_1,b) = (x,b)(c_2,b). \]
		${\cal M}(C; B; P')$ is left cancellative, thus we get $c_1 = c_2$, and that $C$ is left cancellative. \hfill\openbox

		\begin{theorem}Let $A, B, C$ be a right regular triple of semigroups such that $P: B \mapsto A $ is surjective. If $A$ is left reductive, then $C$ is also left reductive.
		\end{theorem}
		
		\textit{Proof.} Assume that $A, B, C$ is a right regular triple of semigroups. Then there are mappings $P: B\mapsto A$ and $P': B\mapsto C$ such that
			\[{\cal M}(A; B; P)/\theta _{{\cal M}(A; B; P)} \cong {\cal M}(C; B; P').\]
			Assume, that $A$ is a left reductive semigroup, and $(a_1,b_1),(a_2,b_2) \in {\cal M}(A; B; P)$ are elements such that
			\[\forall (x,y) \in {\cal M}(A; B; P): \ (x,y)(a_1,b_1) = (x,y)(a_2,b_2).\]
			This means that
			\[ xP(y)a_1 = xP(y)a_2 \quad \textrm{and} \quad b_1 = b_2. \]
			Since $A$ is left reductive, we get that \[\forall y \in B: \ P(y)a_1 = P(y)a_2.\]
			In this case, $P$ is a surjective mapping, hence using again that $A$ is left reductive, we have $a_1 = a_2$. We conclude that $(a_1,b_1) = (a_2,b_2)$, and thus ${\cal M}(A; B; P)$ is left reductive. \\
			\\
			We know, that if $S$ is a left reductive semigroup, then $\theta_S = \iota_S$. This means, that ${\cal M}(A; B; P) \cong {\cal M}(C; B; P')$, hence ${\cal M}(C; B; P')$ is also left reductive.\\
			\\
			Now suppose that $c_1,c_2 \in C$ are such elements, that \[\forall c \in C : \ cc_1 = cc_2. \]
			Take two elements, $(c_1,b), (c_2,b)$ from ${\cal M}(C; B; P')$. For arbitrary $(x,y) \in {\cal M}(C; B; P')$ we have:
			\[ (x,y)(c_1,b) = (xP'(y)c_1,b) = (xP'(y)c_2,b) = (x,y)(c_2,b).\]
			In the second equality, we used the assumption that $\forall c \in C : \ cc_1 = cc_2$. Since ${\cal M}(C; B; P')$ is left reductive, we have $(c_1,b) = (c_2,b)$, and thus $c_1 = c_2$. We conclude that $C$ is left reductive. \hfill\openbox

		\bigskip

		Let $A$ be a semigroup and $B$ be a nonempty set. For a mapping $P$ of $B$ into $A$, let $\alpha _P$ denote the following relation on $A$:
		\[  \alpha_P = \{ (a_1,a_2) \in A \times A : \ (\forall a \in A)(\forall b \in B) \ aP(b)a_1 = aP(b)a_2 \}. \] It is clear that $\alpha_P$ is a right congruence on $A$.  
		
		\begin{remark}\rm 
			It is clear that if $P$ is a mapping of a semigroup $B$ into a semigroup $A$ such that $\alpha _P$ is the identity relation on $A$, then $\theta _{\mathcal{M}(A;B;P)}$ is the identity relation on $\mathcal{M}(A;B;P)$, and hence the triple $A, B, A$ is right regular.
		\end{remark}

\medskip

Let $A, B, C$ be semigroups and $P: B\rightarrow A$, $P':B\rightarrow C$ be arbitrary mappings. We shall say that the triple $A, B, C$ is right regular with respect to the couple $(P, P')$ if $\mathcal{M}(A; B; P)/\theta _{\mathcal{M}(A; B; P)}\cong \mathcal{M}(C; B; P')$.

		\begin{theorem}\label{Th1}
			Let $A$ and $B$ be arbitrary semigroups, and $P$ be a mapping of $B$ into $A$ such that $\alpha_P$ is a congruence on $A$. Then the triple $A, B, A/\alpha _P$ is right regular with respect to $(P, P')$, where
		$P'$ is defined by $P':b \mapsto [P(b)]_{\alpha_P}$ for every $b \in B$. 
		\end{theorem}
		
		\textit{Proof.} Let $\Phi$ be the mapping of the Rees matrix semigroup $M = \mathcal{M}(A;B;P)$ onto the Rees matrix semigroup $\mathcal{M}(A/\alpha_P;B;P')$ defined by
			\[  \Phi: (a,b) \mapsto ([a]_{\alpha_P},b).\]
			For arbitrary elements $(a_1,b_1), (a_2,b_2)$ of $M$, we have
			\[\Phi((a_1,b_1)(a_2,b_2)) = \Phi((a_1P(b_1)a_2,b_2)) = ([a_1P(b_1)a_2]_{\alpha_P},b_2) = \]
			\[ = ([a_1]_{\alpha_P} [P(b_1)]_{\alpha_P}[a_2]_{\alpha_P}, b_2) = ([a_1]_{\alpha_P} P'(b_1)[a_2]_{\alpha_P}, b_2) = \]
			\[ = ([a_1]_{\alpha_P},b_1)([a_2]_{\alpha_P}, b_2) = \Phi((a_1,b_1))\Phi((a_2,b_2)). \]
			Hence, $\Phi$ is a homomorphism. It is clear that $\Phi$ is surjective. We show that the kernel $ker\Phi$ of $\Phi$ is the kernel of the right regular representation of $M$. For elements $(a_1,b_1)$ and $(a_2,b_2)$ of $M$, the equation
			\[ (a,b)(a_1,b_1) = (a,b)(a_2,b_2) \]
			is satisfied for every $a\in A$ and every $b \in B$ if and only if
			\[ (aP(b)a_1, b_1) = (aP(b)a_2,b_2), \]
			that is
			\[ \Phi((a_1,b_1)) = \Phi((a_2,b_2)).\]
			Thus, $ker\Phi = \theta_M$ which proves our theorem. \hfill\openbox

		\vspace{20pt}
		
		A semigroup satisfying the identity $axyb=ayxb$ is called a medial semigroup. It is easy to see that if $A$ is a medial semigroup, then, for an arbitrary semigroup $B$ and an arbirtary mapping of $B$ into $A$, the right congruence $\alpha_P$ is a congruence on $A$. Thus we have the following corollary.
		
		\begin{corollary}
			Let $A$ be a medial semigroup. Then, for an arbitrary semigroup $B$ and an arbitrary mapping $P$ of $B$ into $A$, the triple $A,B,A/\alpha_P$ is right regular, 
		 where $P'$ is defined in Theorem~\ref{Th1}. 
		\end{corollary}

\medskip

If $\varrho$ is an arbitrary congruence on a semigroup $S$, then  $\varrho ^{*}=\{ (a, b)\in S\times S: (\forall s\in S) (sa, sb)\in \varrho\}$ (defined in \cite{Nagy:sg-8}) is also a congruence on $S$ which is called the \emph{right colon congruence of $\varrho$}.

        \begin{remark}\label{Psurjective}\rm
            If $P$ is a mapping  of a nonempty set $B$ onto a semigroup $A$, then $\alpha _P\supseteq \theta ^{*}_A $. If $P$ is surjective, then $\alpha _P=\theta ^{*}_A$.
        \end{remark}
        
\medskip
Remark~\ref{Psurjective} and Theorem~\ref{Th1} imply the following corollary.        
		
		\begin{corollary}
			Let $A$ be an ideal of a semigroup $B$ such that there is a surjective homomorphism $P$ of $B$ onto $A$. Let $P'$ denote the mapping of $B$ onto $A/\theta_A^{\ast}$ defined in the following way: $P':b \mapsto [P(b)]_{\theta_A^{\ast}}$ for every $b \in B$. Then the triple $A, B, A/\theta _A^{\ast}$ is right regular with respect to $(P, P')$.			 
		\end{corollary}
 
\medskip

\medskip

Since the projective homomorphism $P_A: (a, b)\mapsto a$ of the direct product $A\times B$ of semigroups $A$ and $B$ is surjective, Remark~\ref{Psurjective} and Theorem~\ref{Th1} imply the following corollary.
		 
\begin{corollary}
			For arbitrary semigroups $A$ and $B$, the triple $A, A \times B, A/\theta_A^{\ast}$ is right regular with respect to the couple $(P_A, P')$, where $P_A$ denotes the projection homomorphism $P_A : (a,b) \mapsto a$ and $P': A \times B \rightarrow A/\theta_A^{\ast}$ is defined by $P':(a,b) \mapsto [a]_{\theta_A^{\ast}}$.
		\end{corollary}

		\begin{theorem}\label{thmvarphi}
			Let $A$ and $B$ be arbitrary semigroups, and $\varphi$ be a mapping of $A$ into $B$ such that $\alpha _{\varphi}$ is a congruence on $B$. Then the triple $A\times B, A, A/\theta_A^{\ast} \times B/\alpha_{\varphi}$ is right regular with respect to the couple $(P_A, P')$, where $P_A$ is defined by $P_A: a \mapsto (a, \varphi(a))$ and $P'$ is defined by $P':a \mapsto ([a]_{\theta_A^{\ast}},[\varphi(a)]_{\alpha_{\varphi}})$.  
		\end{theorem}
		
		\textit{Proof.} Suppose that $(((a_1,b_1), a_2) , ((a_3,b_3),a_4)) \in \theta_M$, where  \\ $M = \mathcal{M}(A \times B; A; P_A)$. This means that, for every $x, x' \in A$ and $y \in B$,
			\[ ((x,y), x')((a_1,b_1),a_2) = ((x,y),x')((a_3,b_3),a_4)  \iff \]
			\[\iff ((xx'a_1,y\varphi(x')b_1), a_2) = ((xx'a_3, y\varphi(x')b_3), a_4 ).\]
			The equality holds if and only if \[xx'a_1 = xx'a_3,\quad y\varphi(x')b_1 = y\varphi(x')b_3,\quad a_2=a_4,\] that is
\begin{equation}\label{ccc}
(a_1,a_3) \in \theta_A^{\ast},\quad  (b_1, b_3) \in \alpha_{\varphi}, \quad a_2=a_4
\end{equation}

			Let $\Phi$ be the mapping of $\mathcal{M}(A\times B;A;P_A)$ into $\mathcal{M}(A/\theta_A^{\ast} \times B/\alpha_{\varphi};P')$ defined by $\Phi: ((a,b),a') \mapsto (([a]_{\theta_A^{\ast}}, [b]_{\alpha_{\varphi}}), a')$ for every $a,a' \in A$ and every $b \in B$. Since 
			\[ \Phi(((a_1,b_1), a_2)((a_3,b_3),a_4))= \Phi((a_1a_2a_3, b_1 \varphi(a_2)b_3), a_4) =  \]
			\[  = (([a_1a_2a_3]_{\theta_A^{\ast}}, [b_1\varphi(a_2)b_3]_{\alpha_P}), a_4) = (([a_1]_{\theta_A^{\ast}}, [b_1]_{\alpha_{\varphi}}), a_2)(([a_3]_{\theta_A^{\ast}}, [b_3]_{\alpha_{\varphi}}), a_4) =  \]
			\[ = \Phi(((a_1,b_1), a_2)) \Phi(((a_3,b_3),a_4))\] for every $a_1, a_2, a_3, a_4\in A$ and $b_1, b_3\in B$, $\Phi$ is a homomorphism.
			It is clear that $\Phi$ is a surjective. \\
	       Since $(((a_1,b_1),a_2),((a_3,b_3),a_4)) \in ker \Phi$ if and only if all three conditions in (\ref{ccc}) are satisfied, we have $ker \Phi = \theta_M$ and this proves our theorem. \hfill\openbox

		\bigskip

		  If $\varphi: A\mapsto B$ defined in Theorem~\ref{thmvarphi} is surjective, then $\alpha_{\varphi}=\theta _B^{\ast}$ by Remark~\ref{Psurjective}, and thus we have the following corollaries:

\begin{corollary}
			Let $A$ and $B$ be semigroups, and $\varphi$ be a surjective mapping of $A$ onto $B$. Then the triple
			$A\times B, A, A/\theta_A^{\ast} \times B/\theta_B^{\ast}$ is right regular with respect to the couple $(P_A, P')$, where $P_A$ is defined by $P_A: a \mapsto (a, \varphi(a))$ and $P'$ is defined by $P':a \mapsto ([a]_{\theta_A^{\ast}},[\varphi(a)]_{\theta_B^{\ast}})$.
		\end{corollary}
		
		\begin{corollary}
			Let A be a semigroup, and $B$ be a retract ideal of $A$. Let $\varphi$ be a retract homomorphism of $A$ onto $B$. Then the triple
			$A\times B, A, A/\theta_A^{\ast} \times B/\theta_B^{\ast}$ is right regular with respect to the couple $(P_A, P')$, where $P_A$ is defined by $P_A: a \mapsto (a, \varphi(a))$ and $P'$ is defined by $P':a \mapsto ([a]_{\theta_A^{\ast}},[\varphi(a)]_{\theta_B^{\ast}})$.
		\end{corollary}
		
		\medskip
		
		If $B$ is an ideal of a semigroup $A$ such that $B$ is a group, then $\varphi _B: A\rightarrow B$ defined by $\varphi_B (a)=ae$ ($a\in A$) is a retract homomorphism of $A$ onto $B$, where $e$ denotes the identity element of the group $B$.  
		
		\begin{corollary}
			Let $A$ be a semigroup and $B$ be an ideal of $A$ such that $B$ is a group. Then the triple $A\times B, A, A/\theta_A^{\ast} \times B$ is right regular with respect to the couple $(P_A, P')$, where $P_A$ is defined by $P_A: a \mapsto (a, \varphi _B(a))$ and $P'$ is defined by $P':a \mapsto ([a]_{\theta_A^{\ast}},\varphi _B(a))$; here $\varphi _B$ denotes the above surjective homomorphism of $A$ onto $B$.
		\end{corollary}

     \end{document}